\documentclass[english]{article}
\usepackage[T1]{fontenc}
\usepackage[latin9]{inputenc}
\usepackage{amsmath}
\usepackage{amssymb}
\usepackage{stmaryrd}
\usepackage{graphicx}

\makeatletter

\providecommand{\tabularnewline}{\\}

\makeatother

\usepackage{babel}
\begin{document}
\title{Approximation by Power Series of Functions}
\author{Andrej Liptaj\thanks{andrej.liptaj@savba.sk, ORC iD 0000-0001-5898-6608}\\
{\small{}Institute of Physics, Bratislava, Slovak Academy of Sciences}}
\maketitle
\begin{abstract}
Derivative-matching approximations are constructed as power series
built from functions. The method assumes the knowledge of special
values of the Bell polynomials of the second kind, for which we refer
to the literature. The presented ideas may have applications in numerical
mathematics.
\end{abstract}

\section*{Introduction}

Given a function $f$ and a point of expansion $x_{0}$, it is customary
to say that the Taylor polynomial (TP) of degree one, two, three,...
is the best linear, quadratic, cubic,... approximation of $f$ at
$x_{0}$. In this sense we present here several new approximations
$\mathcal{A}_{i}^{f}$ of $f$ such that 
\begin{equation}
\frac{d^{n}}{dx^{n}}f\left(x\right)|_{x=0}=\frac{d^{n}}{dx^{n}}\mathcal{A}_{i}^{f}\left(x\right)|_{x=0},\quad n\in\mathbb{N}_{0},\label{eq:derMatching}
\end{equation}
where, without loss of generality, we assume that the expansion is
done at $x_{0}=0$ (shift to an arbitrary point $x_{0}$ is achieved
by shifting the argument). We denote the equality (\ref{eq:derMatching})
by $f\approx\mathcal{A}_{i}^{f}$.

\section{Power series built from functions}

We build $\mathcal{A}_{i}^{f}$ as a power series of some properly
chosen function $g$ following the construction from Sec. 4.1.2 of
\cite{jaArxiv}. We propose 
\begin{equation}
\mathcal{A}_{i}^{f}\left(x\right)=\sum_{n=0}^{\infty}a_{n}\left[g\left(x\right)\right]^{n}\approx f(x)\text{ with }\;g\left(0\right)=0\text{ and }g'\left(0\right)\neq0.\label{eq:xLikeExpansion}
\end{equation}
The existence of a non-zero derivative at zero implies $g$ can be
inverted on some neighborhood of zero $x\equiv g^{-1}\left(y\right)$.
We have
\begin{equation}
f\left[g^{-1}\left(y\right)\right]\approx\sum_{n=0}^{\infty}a_{n}y^{n},\label{eq:Substitucia}
\end{equation}
i.e. the expansion coefficients $a_{n}$ are given by the power expansion
coefficients of $f\left(g^{-1}\right)$
\begin{equation}
a_{n}=\frac{1}{n!}\frac{d^{n}}{dx^{n}}f\left(g^{-1}\left(x\right)\right)|_{x=0}.\label{eq:LukasovVzorec}
\end{equation}
This can be written in terms of the Faà di Bruno's formula, where
the Bell polynomials of the second kind $B_{n,k}$ appear
\begin{equation}
a_{n}=\frac{1}{n!}\sum_{k=0}^{n}d_{k}^{f}B_{n,k}(d_{1}^{g^{-1}},d_{2}^{g^{-1}},\ldots,d_{n-k+1}^{g^{-1}});\quad d_{n}^{h}\equiv\frac{d^{n}}{dx^{n}}h\left(x\right)|_{x=0}.\label{eq:koeficentyVzorec}
\end{equation}
In \cite{jaArxiv} only few expansions were presented, here we systematically
review the existing formulas for special values of the Bell polynomials
\cite{bellSpecVal1,bellSpecVal2,bellSpecVal3} and propose a larger
number of them\footnote{Included are also those from \cite{jaArxiv}, so as to provide a complete
list of approximations of this kind.}.

To keep the text brief, we organize our results as a list where only
the necessary information is summarized. We define
\begin{align*}
(-1)!! & =1,\qquad0^{0}=1,\qquad\left\langle \alpha\right\rangle _{n}=\prod_{k=0}^{n-1}\left(\alpha-k\right)\text{ //falling factorial,}\\
W\left(x\right) & \rightarrow\text{principal branch of the Lambert W function,}\\
\left\llbracket \begin{array}{c}
n\\
m
\end{array}\right\rrbracket  & =\frac{1}{m!}\sum_{k=0}^{m}\left(-1\right)^{k}\binom{m}{k}\left(m-k\right)^{n},\\
 & \text{ (Stirling numbers of the second kind)}\\
\left[\begin{array}{c}
n\\
m
\end{array}\right] & =\sum_{j=0}^{n-m}\left(-1\right)^{j}\binom{n-1+j}{n-m+j}\binom{2n-m}{n-m-j}\left\llbracket \begin{array}{c}
n-m+j\\
j
\end{array}\right\rrbracket .\\
 & \text{ (Stirling numbers of the first kind)}
\end{align*}
When needed, we extend the definition of $g$ (or $g^{-1}$) to zero
by its limit value
\[
g(0)=\lim_{x\rightarrow0^{(\pm)}}g\left(x\right),
\]
and note it with $\doteq$. The exact version of the limit (left,
right, both sides) depends on the context.

\section{List of expansions}

The expansion is for all cases constructed as
\begin{equation}
\mathcal{A}_{i}^{f}\left(x\right)=f(0)+\sum_{n=1}^{N}\frac{1}{n!}\left[\sum_{k=1}^{n}d_{k}^{f}B_{n,k}(d_{1}^{g^{-1}},d_{2}^{g^{-1}},\ldots,d_{n-k+1}^{g^{-1}})\right]\left[g\left(x\right)\right]^{n},\label{eq:mainSeries}
\end{equation}
where we isolate the constant term so as to avoid ambiguities for
$n=0$ (such as $0^{0}$) in the formulas which follow. We separate
cases where an explicit formula for $g$ is found and those where
it is not. In the first scenario we present also the formula for the
Belle polynomial values\footnote{We want to provide the full information needed for an eventual implementation,
so that the reader does not need to look into the literature we cite.}, in the second situation we do this only for short formulas, for
the long ones we cite the literature. The displayed constants directly
appearing as arguments of the Belle polynomials $B_{n,k}(c_{1},c_{2},c_{3},\ldots)$
give the information about the derivatives of $g^{-1}$ at zero for
the case in question, i.e. $c_{i}=d_{i}^{g^{-1}}$.

\subsection{Formulas with explicit expression for $g$\label{subsec:ExplicitG}}
\begin{itemize}
\item Logarithm-based expansion ($\mathcal{A}_{1}^{f}$)
\begin{align}
g\left(x\right) & =\ln\left(x+1\right);\quad g^{-1}\left(x\right)=\exp\left(x\right)-1,\label{eq:app1}
\end{align}
\[
B_{n,k}(1,1,1,\ldots)=\left\llbracket \begin{array}{c}
n\\
k
\end{array}\right\rrbracket .
\]
\item Exponential-based expansion ($\mathcal{A}_{2}^{f}$)
\begin{align}
g\left(x\right) & =1-e^{-x};\quad g^{-1}\left(x\right)=-\ln\left(1-x\right),\label{eq:app2}
\end{align}
\[
B_{n,k}(0!,1!,2!,\ldots)=\left(-1\right)^{n-k}\left[\begin{array}{c}
n\\
k
\end{array}\right].
\]
\item Expansion with inverse hyperbolic sine ($\mathcal{A}_{3}^{f}$)
\begin{align}
g\left(x\right) & =\text{asinh}(x);\quad g^{-1}\left(x\right)=\sinh(x),\label{eq:app3}
\end{align}
\[
B_{n,k}(1,0,1,0,1\ldots)=\frac{1}{2^{k}k!}\sum_{l=0}^{k}\left(-1\right)^{l}\binom{k}{l}\left(k-2l\right)^{n}.
\]
\item Arcus-sine-based expansion ($\mathcal{A}_{4}^{f}$)
\begin{align}
g\left(x\right) & =\arcsin(x);\quad g^{-1}\left(x\right)=\sin(x),\label{eq:app4}
\end{align}
\[
B_{n,k}(1,0,-1,0,1\ldots)=\frac{\left(-1\right)^{k}}{2^{k}k!}\cos\left[\frac{\left(n-k\right)\pi}{2}\right]\sum_{q=0}^{k}\left(-1\right)^{q}\binom{k}{q}\left(2q-k\right)^{n}.
\]
\item Expansion in powers of $\sqrt[\alpha]{x+1}-1$ ($\mathcal{A}_{5}^{f}$)
\begin{align}
g\left(x\right) & =\sqrt[\alpha]{x+1}-1;\quad g^{-1}\left(x\right)=(1+x)^{\alpha}-1;\quad\alpha\in\mathbb{R}\backslash\left\{ 0\right\} ,\label{eq:app5}
\end{align}
\[
B_{n,k}\left(\left\langle \alpha\right\rangle _{1},\left\langle \alpha\right\rangle _{2},\left\langle \alpha\right\rangle _{3},\ldots\right)=\frac{\left(-1\right)^{k}}{k!}\sum_{l=0}^{k}\left(-1\right)^{l}\binom{k}{l}\left\langle \alpha l\right\rangle _{n}.
\]
Notable spacial cases (polynomial and rational) happen for $\alpha=\pm1/n$,
$n\in\mathbb{N}$. For $\alpha=1$ the TP is constructed.
\item Square-root-based expansion ($\mathcal{A}_{6}^{f}$)
\begin{align}
g\left(x\right) & =\sqrt{2x+w^{2}}-w;\quad g^{-1}\left(x\right)=\frac{1}{2}x^{2}+wx;\quad w\in\mathbb{R}\backslash\left\{ 0\right\} ,\label{eq:app6}
\end{align}
\[
B_{n,k}\left(w,1,0,0,0,\ldots\right)=\frac{1}{2^{n-k}}\frac{n!}{k!}\binom{k}{n-k}w^{2k-n}.
\]
\item Polynomial expansion ($\mathcal{A}_{7}^{f}$)
\begin{align}
g\left(x\right) & =\frac{x^{2}+2\sqrt{\alpha}x}{\beta};\quad g^{-1}\left(x\right)=\sqrt{\alpha+\beta x}-\sqrt{\alpha};\quad\alpha,\beta\in\mathbb{R}\backslash\left\{ 0\right\} ,\label{eq:app7}
\end{align}
\[
B_{n,k}\left(d_{1}^{g^{-1}},d_{2}^{g^{-1}},\ldots\right)=\left(-1\right)^{n+k}\frac{\left[2\left(n-k\right)-1\right]!!}{\alpha^{n-k/2}}\left(\frac{\beta}{2}\right)^{n}\binom{2n-k-1}{2\left(n-k\right)},
\]
where
\[
d_{n}^{g^{-1}}=\alpha^{\frac{1}{2}-n}\beta^{n}\prod_{k=1}^{n}\left(k+\frac{1}{2}-n\right).
\]
\item Expansion with the square root in the denominator ($\mathcal{A}_{8}^{f}$)
\begin{align}
g\left(x\right) & =1-\frac{1}{\sqrt{x+1}};\quad g^{-1}\left(x\right)=\frac{1}{\left(x-1\right)^{2}}-1,\label{eq:app8}
\end{align}
\[
B_{n,k}(2!,3!,4!,\ldots)=\frac{n!}{k!}\sum_{l=0}^{k}\left(-1\right)^{k-l}\binom{k}{l}\binom{n+2l-1}{n}.
\]
\item Expansion with fraction including square root ($\mathcal{A}_{9}^{f}$)
\begin{align}
g\left(x\right) & \doteq\frac{-1+\sqrt{4x^{2}+1}}{2x};\quad g^{-1}\left(x\right)=\frac{x}{1-x^{2}},\label{eq:app9}
\end{align}
\[
B_{n,k}(1!,0,3!,0,5!,0\ldots)=\frac{1+\left(-1\right)^{n+k}}{2}\frac{n!}{k!}\binom{\frac{n+k}{2}-1}{k-1}.
\]
\item Expansion with the Lambert function ($\mathcal{A}_{10}^{f}$)
\begin{align}
g\left(x\right) & =W\left[e^{w-1}\left(w+x-1\right)\right]+1-w,\quad w\in\mathbb{R}\backslash\left\{ 0\right\} ,\label{eq:app10}\\
g^{-1}\left(x\right) & =\left(w+x-1\right)e^{x}+1-w,\nonumber 
\end{align}
\begin{align*}
B_{n,k}(w,w+1,w+2,\ldots) & =\\
\qquad=k^{n-k}\binom{n}{k}\sum_{l=0}^{k}\binom{k}{l} & \left[\sum_{q=0}^{n-k}\frac{\left(-1\right)^{q}}{k^{q}}\binom{n-k}{q}\frac{\left\llbracket \begin{array}{c}
l+q\\
l
\end{array}\right\rrbracket }{\binom{l+q}{l}}\right]\left(w-1\right)^{l}.
\end{align*}
\item Second expansion with the Lambert function ($\mathcal{A}_{11}^{f}$)
\begin{align}
g\left(x\right) & =\frac{W\left[-e^{-(x+1)}(x+1)\right]}{x+1}+1;\quad g^{-1}\left(x\right)\doteq-\frac{\ln\left(1-x\right)}{x}-1,\label{eq:app11}
\end{align}
\[
B_{n,k}\left(\frac{1!}{2},\frac{2!}{3},\frac{3!}{4},\ldots\right)=\frac{(-1)^{n-k}}{k!}\sum_{m=0}^{k}\left(-1\right)^{m}\binom{k}{m}\frac{\left[\begin{array}{c}
n+m\\
m
\end{array}\right]}{\binom{n+m}{m}}.
\]
As readily seen form the argument of the function $W$(which is defined
from $-1/e$ to $\infty$), this approximation is valid in the right
neighborhood of zero.
\item Third expansion with the Lambert function ($\mathcal{A}_{12}^{f}$)
\begin{align}
g\left(x\right) & =-\frac{W\left(-\frac{\exp\left(-\frac{1}{1+x}\right)}{1+x}\right)+xW\left(-\frac{\exp\left(-\frac{1}{1+x}\right)}{1+x}\right)+1}{1+x},\label{eq:app12}\\
g^{-1}\left(x\right) & \doteq\frac{e^{x}-1}{x}-1,\nonumber 
\end{align}
\[
B_{n,k}\left(\frac{1}{2},\frac{1}{3},\frac{1}{4},\ldots\right)=\frac{n!}{\left(n+k\right)!}\sum_{l=0}^{k}\left(-1\right)^{k-l}\binom{n+k}{k-l}\left\llbracket \begin{array}{c}
n+l\\
l
\end{array}\right\rrbracket .
\]
As readily seen form the argument of the function $W$, this approximation
is valid in the left neighborhood of zero.
\item Powers of sine ($\mathcal{A}_{13}^{f}$)
\begin{align}
g\left(x\right) & =\sin\left(x\right),\label{eq:app13}\\
g^{-1}\left(x\right) & =\arcsin\left(x\right),\nonumber 
\end{align}
\begin{align*}
B_{n,k}\left(1,0,1,0,9,0,225,0,\ldots,\left[\left(n-k-3\right)!!\right]^{2},0,\left[\left(n-k-1\right)!!\right]^{2}\right)=\\
=\delta_{\left(n-k\right)\%2,0}\left(-1\right)^{\frac{n-k}{2}}2^{n-k}\sum_{l=0}^{n-k}\binom{k+l-1}{k-1}\left[\begin{array}{c}
n-1\\
k+l-1
\end{array}\right]\left(\frac{n-2}{2}\right)^{l},
\end{align*}
where $\delta$ is the Kronecker delta and $\%$ is the modulo operation.
This expansion has large similarities with \cite{Butzer} and represents
Fourier series whose standard form can be get by applying trigonometric
power formulas to $\left[\sin\left(x\right)\right]^{n}$ terms.
\end{itemize}

\subsection{Formulas without explicit expression for $g$}

With the function $g^{-1}$ known, one can use numerical or approximation
methods to get $g$ in the proximity of zero.
\begin{itemize}
\item Case one\\
\begin{equation}
g^{-1}\left(x\right)=(w-1+e^{x})x;\quad w\neq0,\label{eq:app14}
\end{equation}
\[
B_{n,k}(w,2,3,4,\ldots)=\binom{n}{k}\sum_{r=0}^{k}\binom{k}{r}\left(k-r\right)^{n-k}\left(w-1\right)^{r}.
\]
\item Case two\\
\begin{equation}
g^{-1}(x)=e^{x}(x-2)-x+2,\label{eq:app15}
\end{equation}
\begin{align*}
B_{n,k}(-2,0,1,2,3,\ldots) & =\sum_{r=0}^{n}r!\binom{n}{r}\binom{k}{r}\left(-2\right)^{k-r}\left\llbracket \begin{array}{c}
n-r\\
k
\end{array}\right\rrbracket .
\end{align*}
\item Case three\\
\begin{equation}
g^{-1}(x)=\left(2e^{x}-x^{2}-2x-2\right)/\left(2x^{2}\right).\label{eq:app16}
\end{equation}
The formula for $B_{n,k}(\frac{1}{2.3},\frac{1}{3.4},\ldots)$ is
shown in Eq. (2.1) of \cite{bellSpecVal1}.
\item Case four\\
\begin{equation}
g^{-1}(x)=\left(6xe^{x}-12e^{x}-x^{3}+6x+12\right)/\left(6x^{3}\right).\label{eq:app17}
\end{equation}
The formula for $B_{n,k}(\frac{1}{3.4},\frac{1}{4.5},\ldots)$ is
shown in Theorem 2.7 of \cite{bellSpecVal1}.
\item Case five\\
\begin{equation}
g^{-1}(x)=\alpha+(\alpha+a_{1}-1)x+\frac{1}{2}(\alpha+a_{2}-2)x^{2}+(x-\alpha)e^{x};\quad a_{1}\neq0.\label{eq:app18}
\end{equation}
The formula for $B_{n,k}(a_{1},a_{2},3-\alpha,4-\alpha,5-\alpha,\ldots)$
is shown in Eq. (3.1) of \cite{bellSpecVal1}. The function $g$ can
be expressed in terms of the Lambert $W$ for $a_{1}=1-\alpha$ and
$a_{2}=2-\alpha$, which however corresponds to Eq. (\ref{eq:app10})
from the previous section.
\item Case six\\
\begin{equation}
g^{-1}(x)\doteq-\frac{\left[\arccos\left(x+1\right)\right]^{2}}{2x}-1.\label{eq:app19}
\end{equation}
The formula for $B_{n,k}\left(-\frac{2}{12},\frac{4}{45},-\frac{6}{70},\ldots,2\frac{(2n-2k+2)!!}{(2n-2k+4)!}Q(2,2n-2k+2)\right)$
together with the definition of $Q$ is shown in Eqs. (5.1) and (2.3)
of \cite{bellSpecVal2}.
\end{itemize}

\section{Discussion and remarks}

\subsection*{Plots}

\begin{figure}
\begin{centering}
\includegraphics[width=0.5\textwidth]{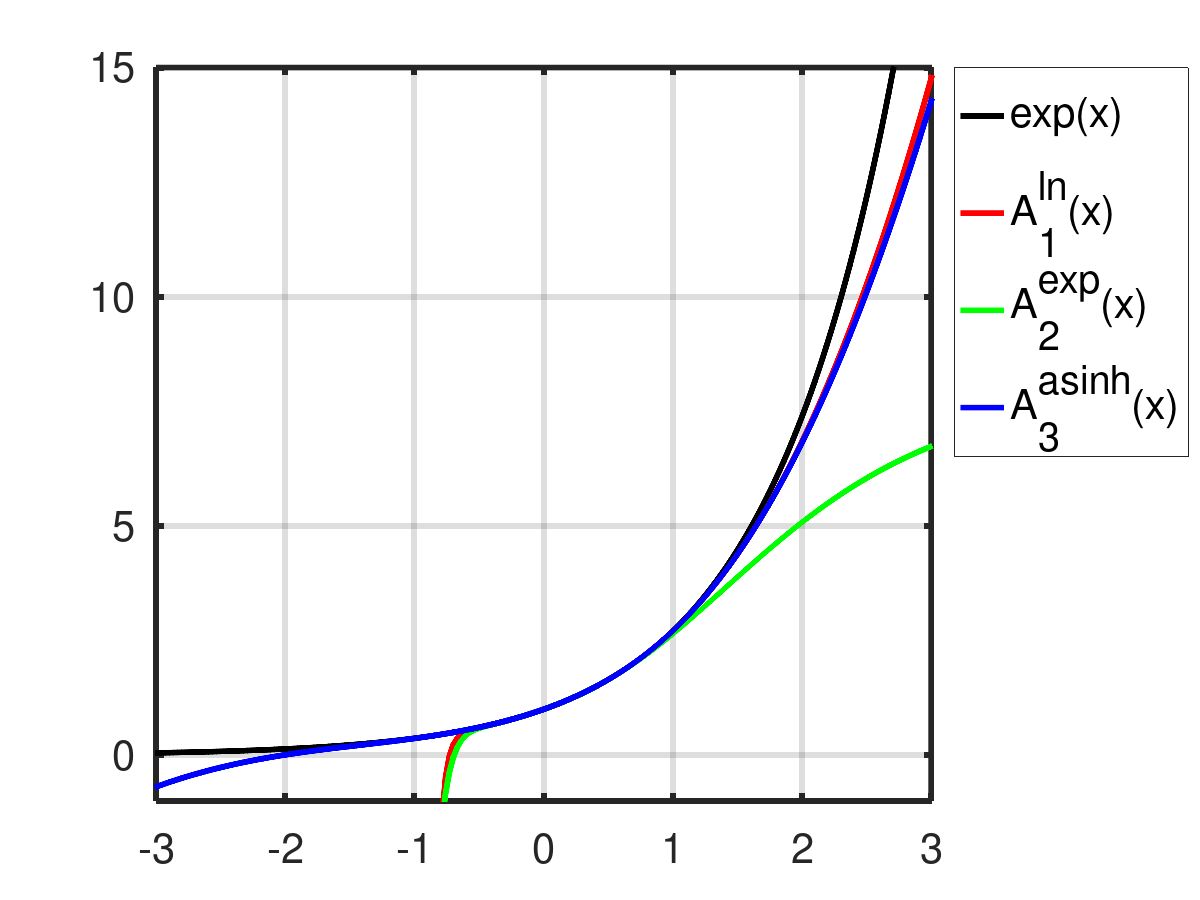}\includegraphics[width=0.5\textwidth]{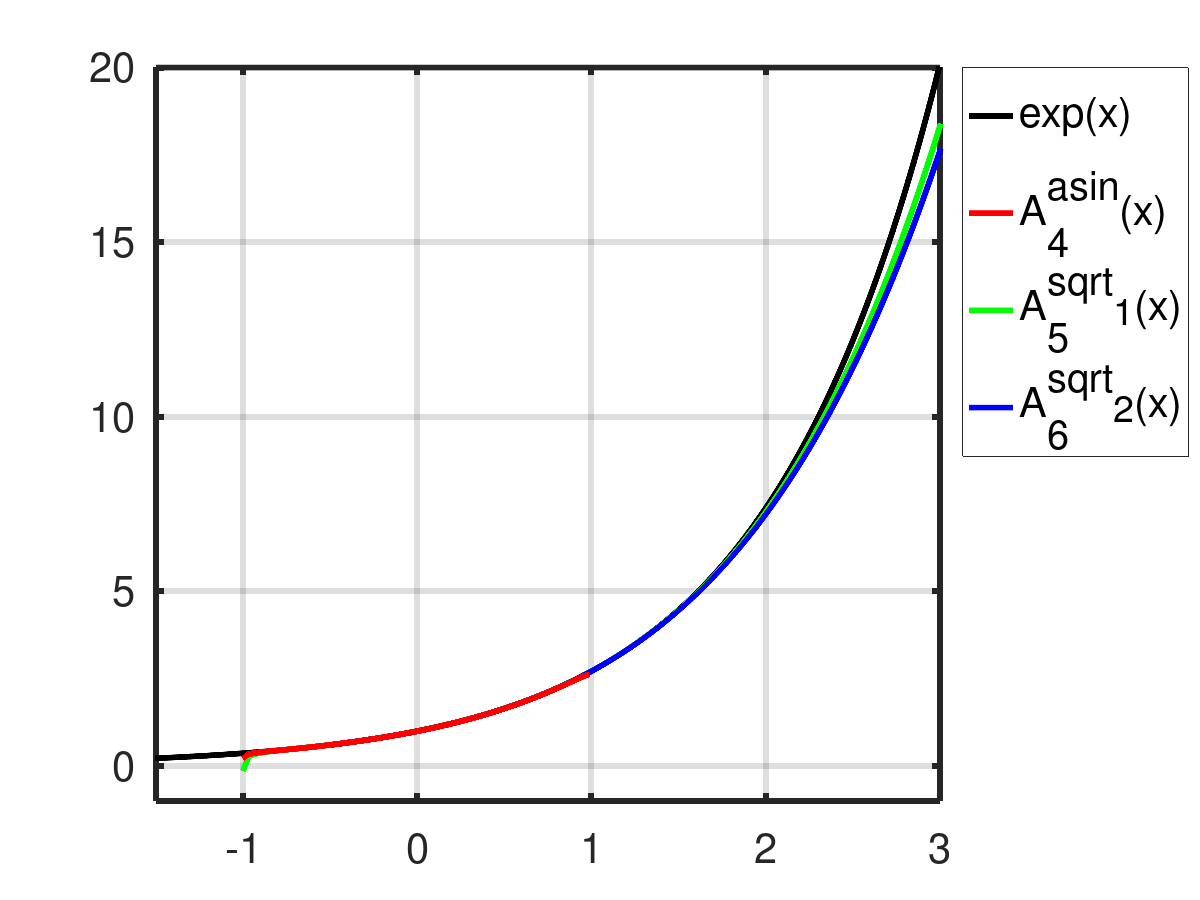}
\par\end{centering}
\begin{centering}
\includegraphics[width=0.5\textwidth]{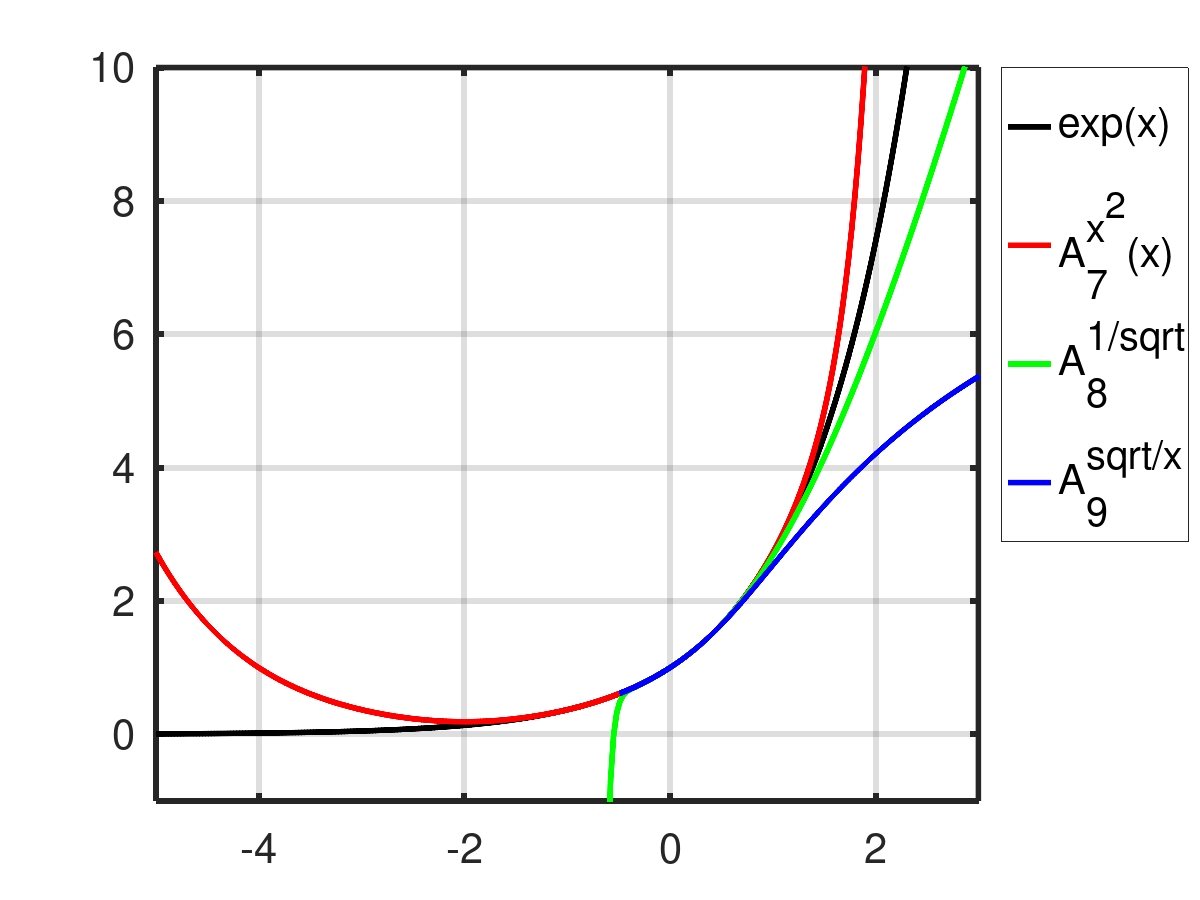}\includegraphics[width=0.5\textwidth]{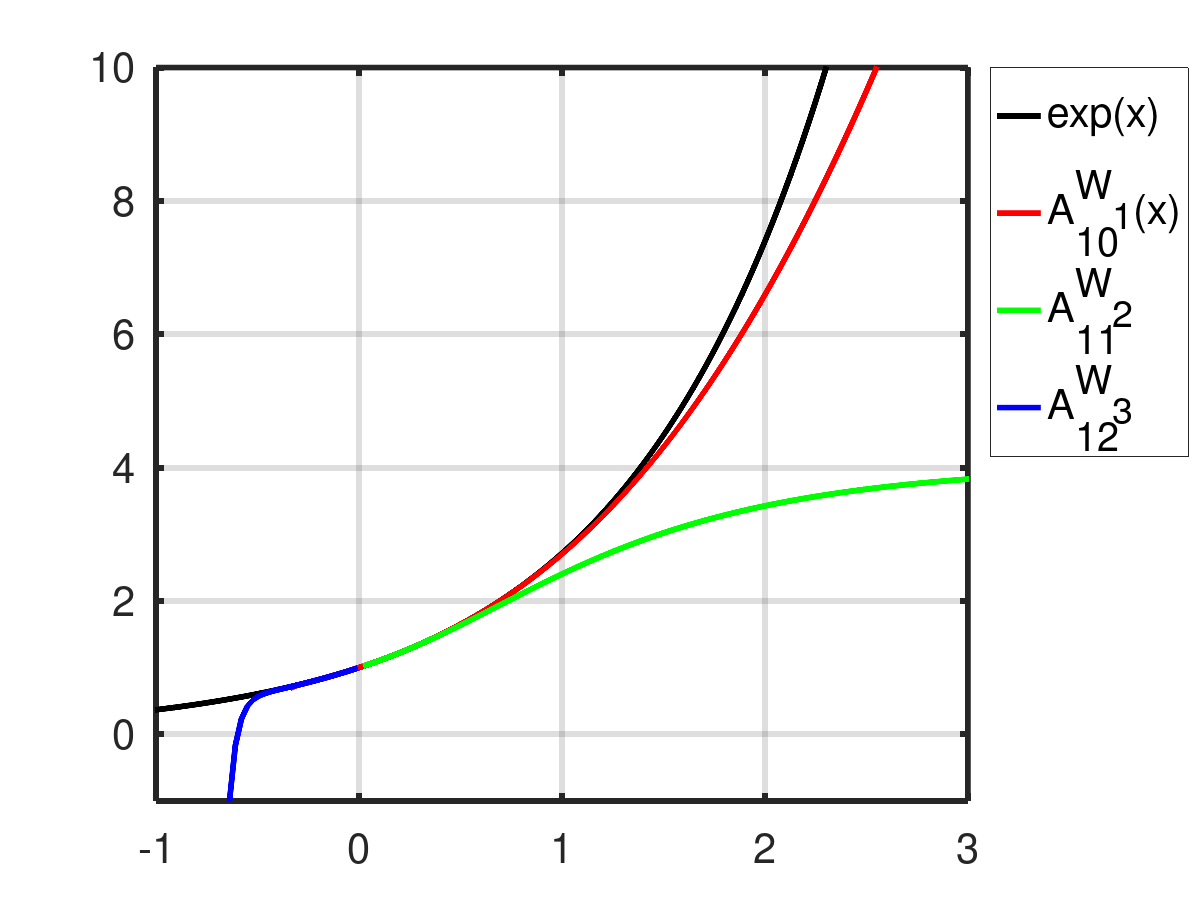}
\par\end{centering}
\begin{centering}
\includegraphics[width=0.5\textwidth]{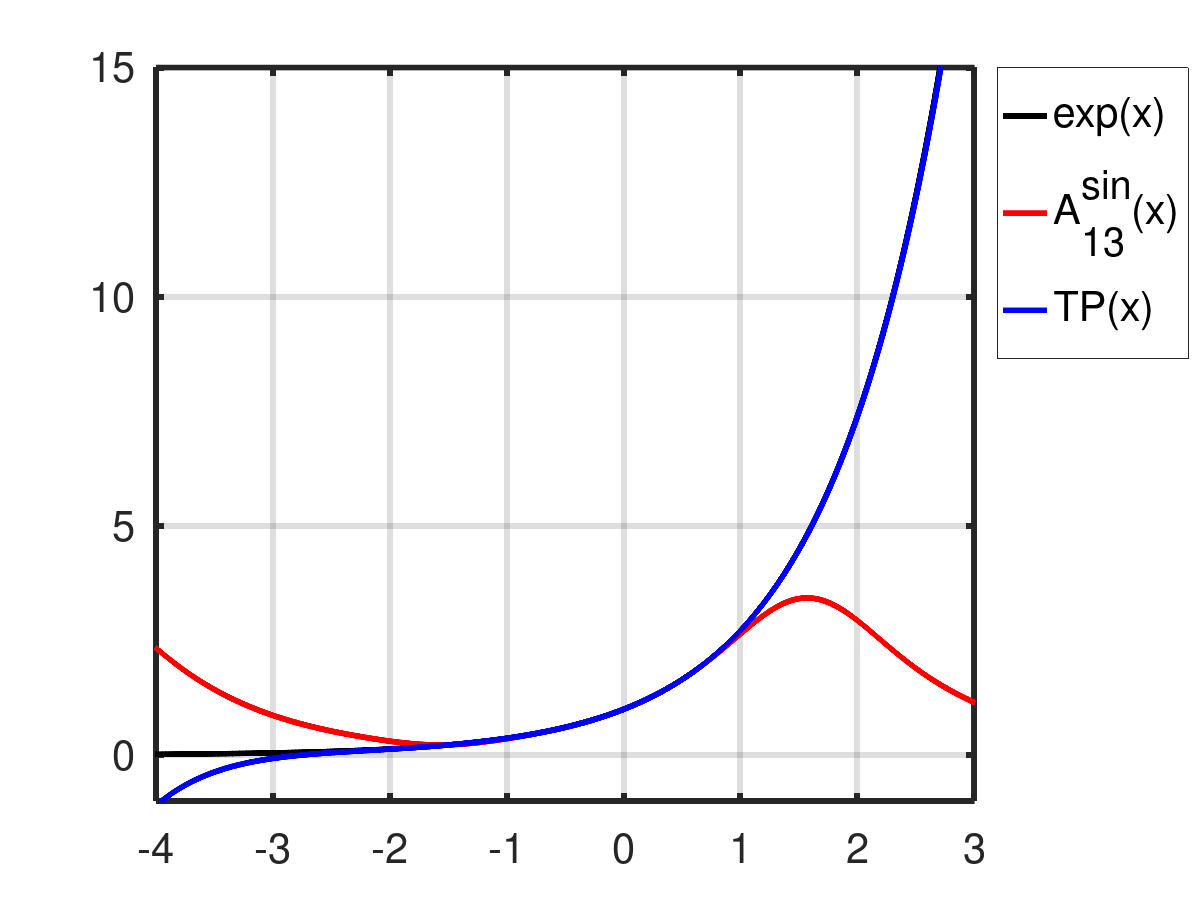}
\par\end{centering}
\caption{The $\exp\left(x\right)$ function approximated by series (\ref{eq:mainSeries})
with 8 terms and with various $g$ from Sec. \ref{subsec:ExplicitG}.}
\label{FigExp}
\end{figure}

\begin{figure}
\begin{centering}
\includegraphics[width=0.5\textwidth]{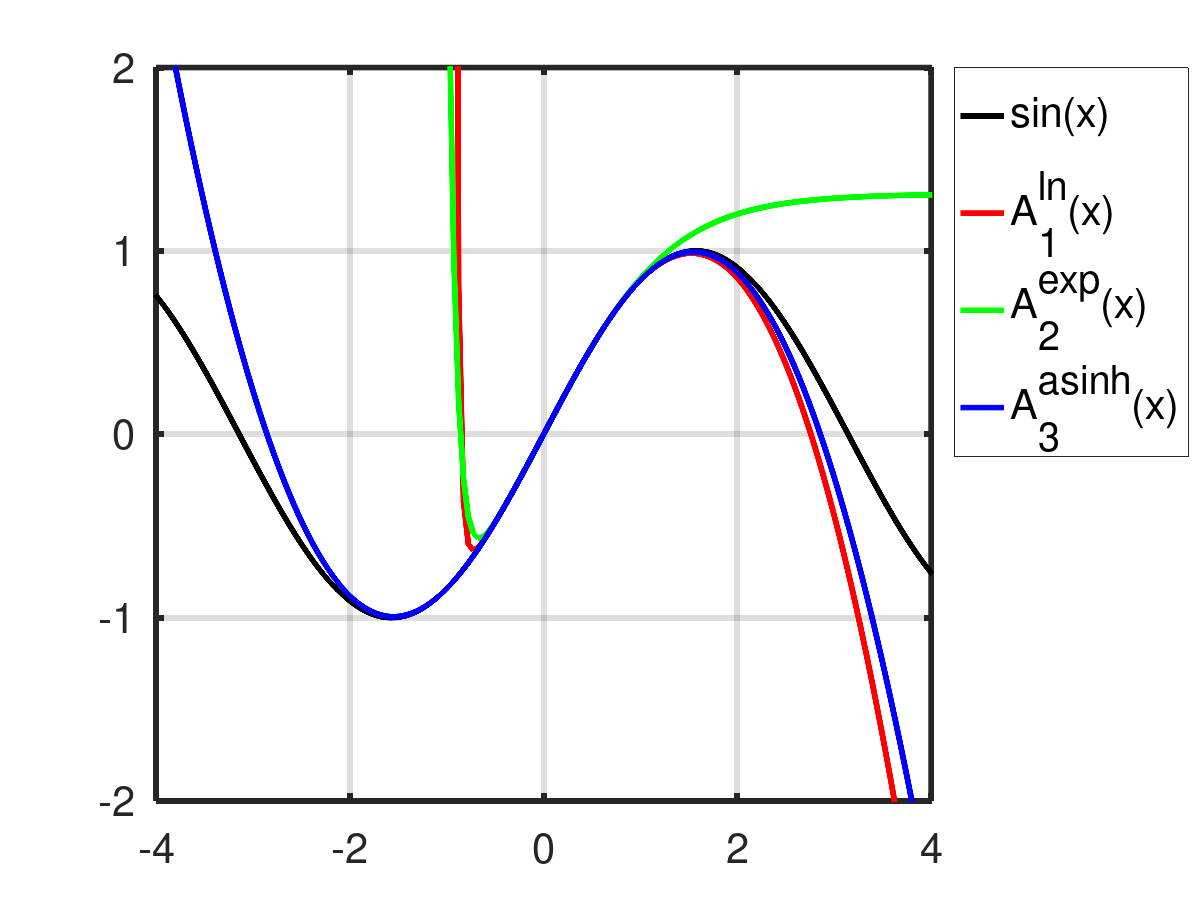}\includegraphics[width=0.5\textwidth]{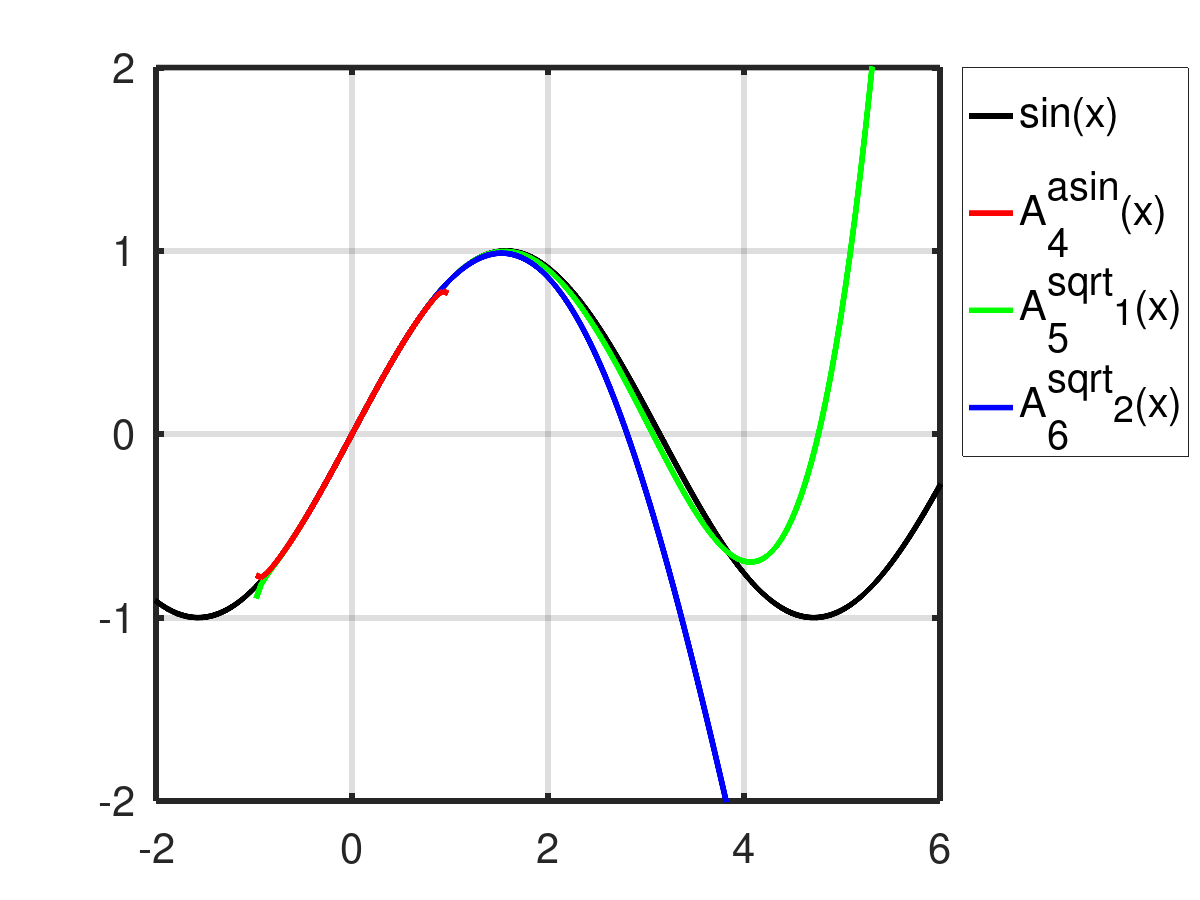}
\par\end{centering}
\begin{centering}
\includegraphics[width=0.5\textwidth]{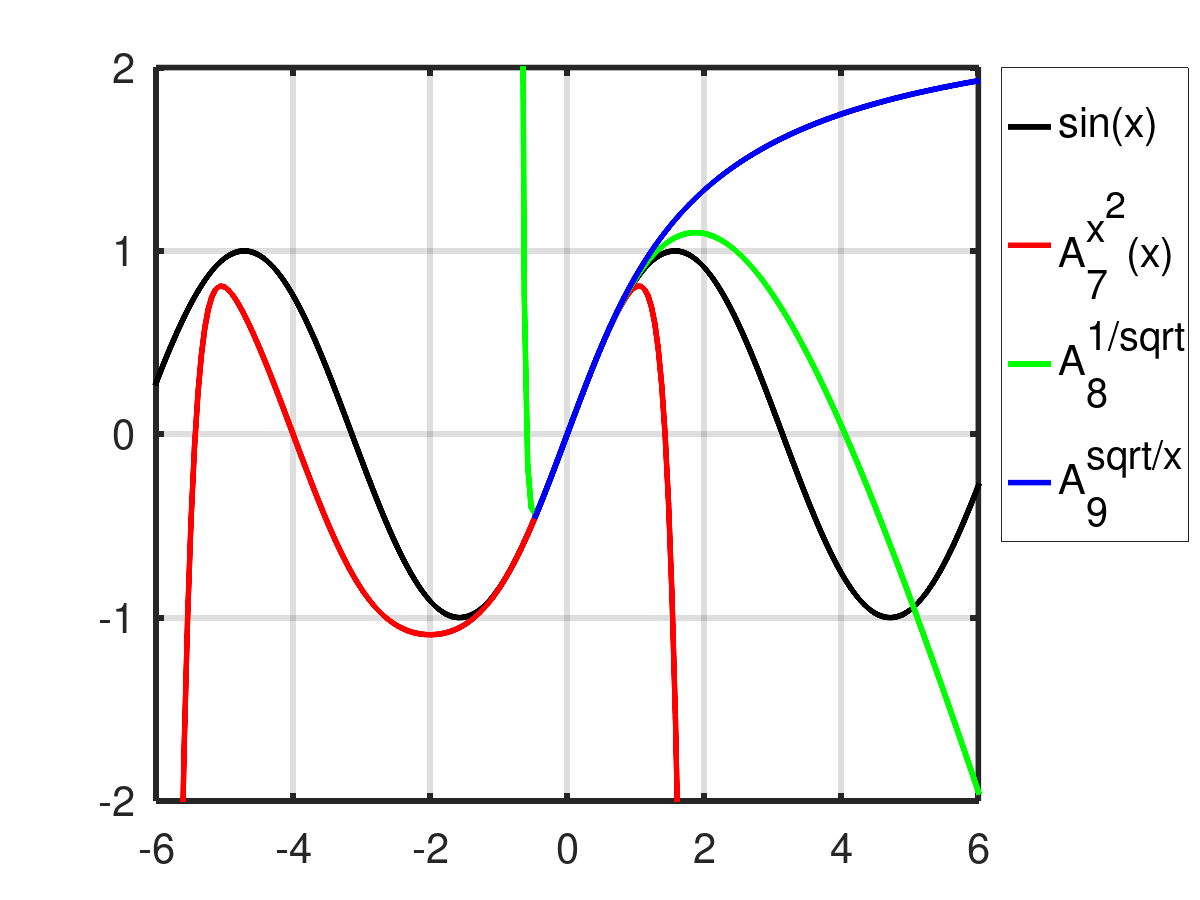}\includegraphics[width=0.5\textwidth]{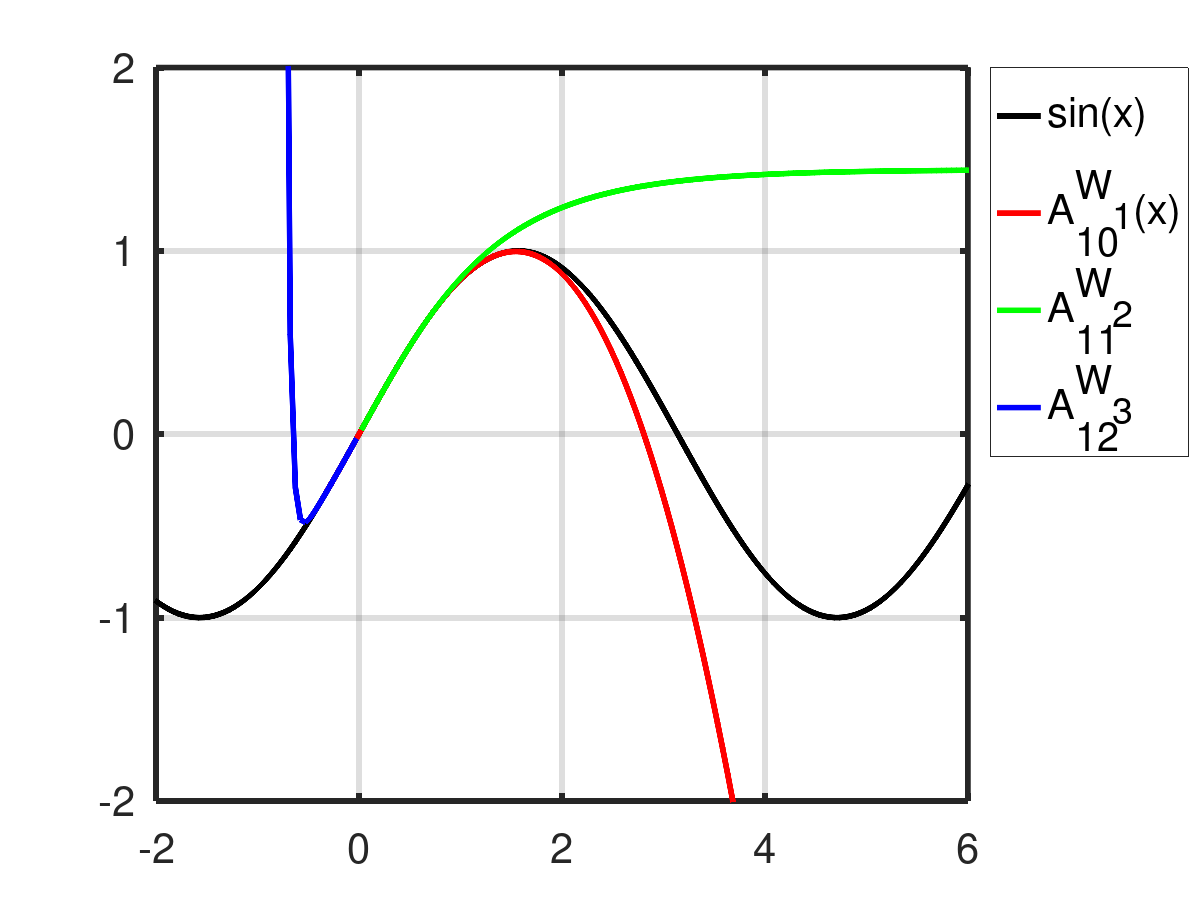}
\par\end{centering}
\begin{centering}
\includegraphics[width=0.5\textwidth]{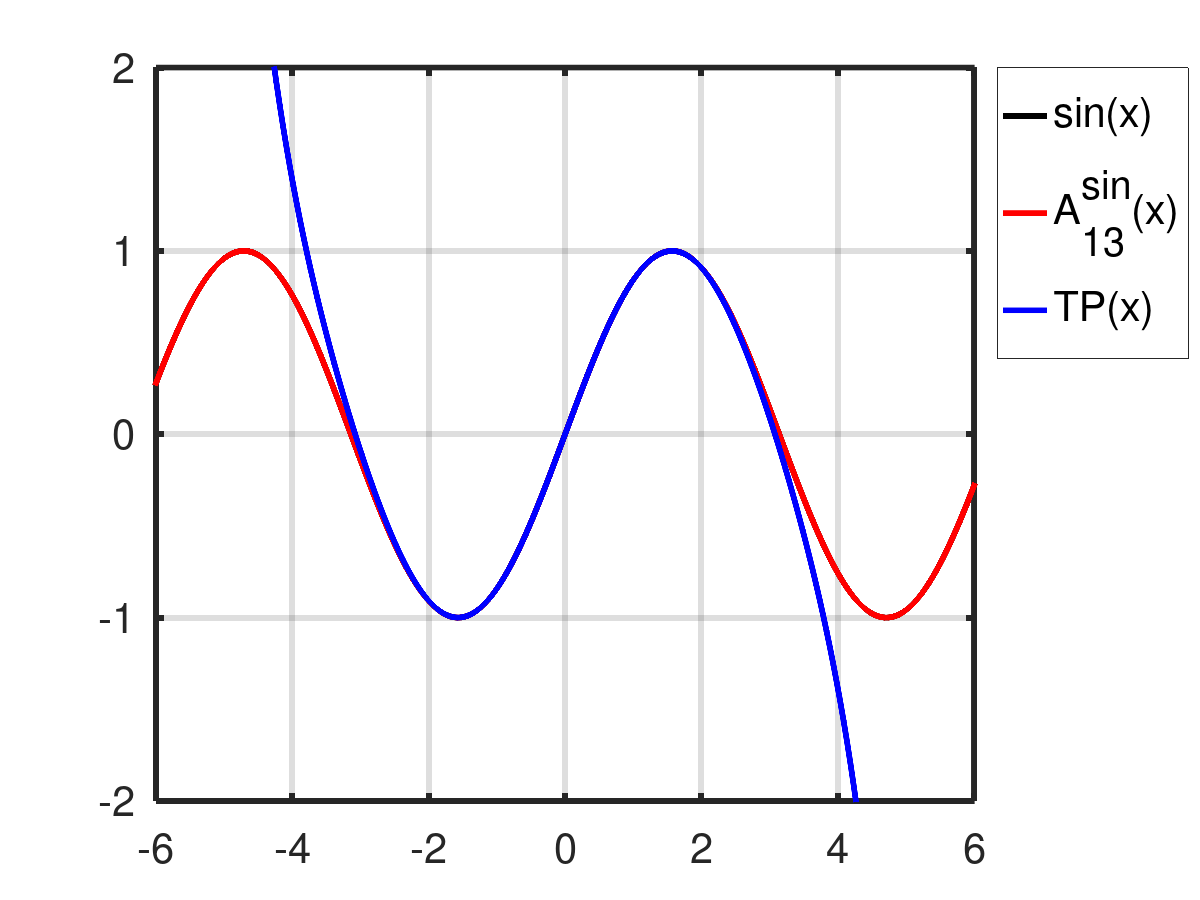}
\par\end{centering}
\caption{The $\sin\left(x\right)$ function approximated by series (\ref{eq:mainSeries})
with 8 terms and with various $g$ from Sec. \ref{subsec:ExplicitG}. }
\label{FigSin}
\end{figure}

\begin{figure}
\begin{centering}
\includegraphics[width=0.5\textwidth]{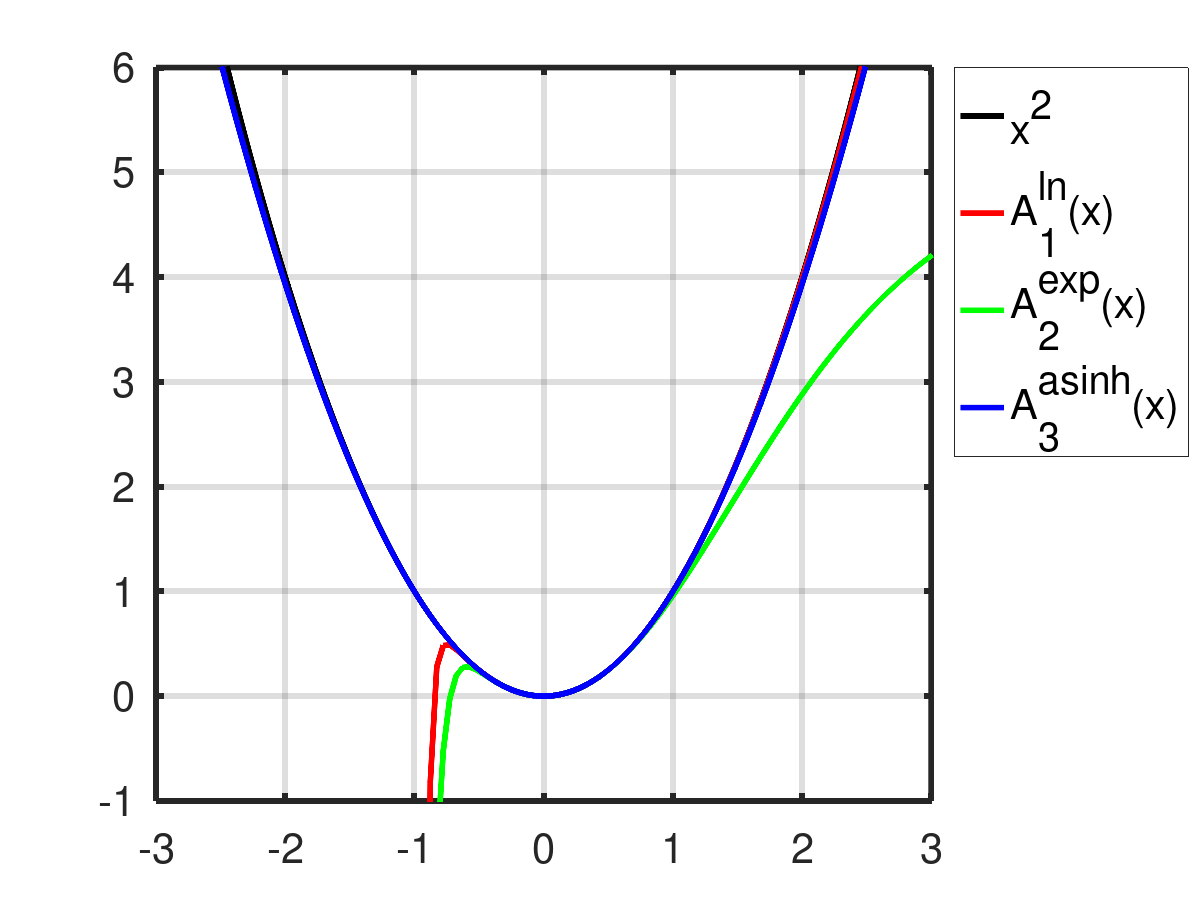}\includegraphics[width=0.5\textwidth]{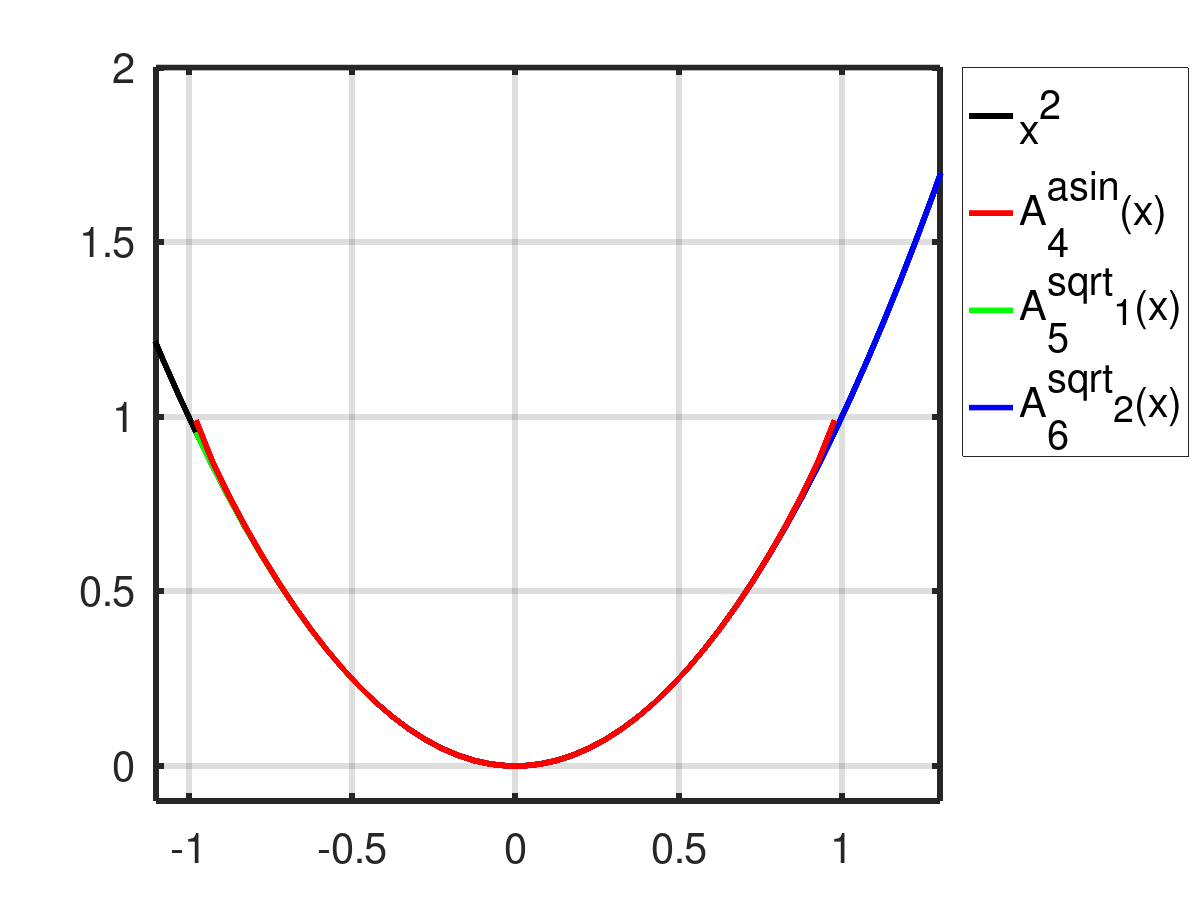}
\par\end{centering}
\begin{centering}
\includegraphics[width=0.5\textwidth]{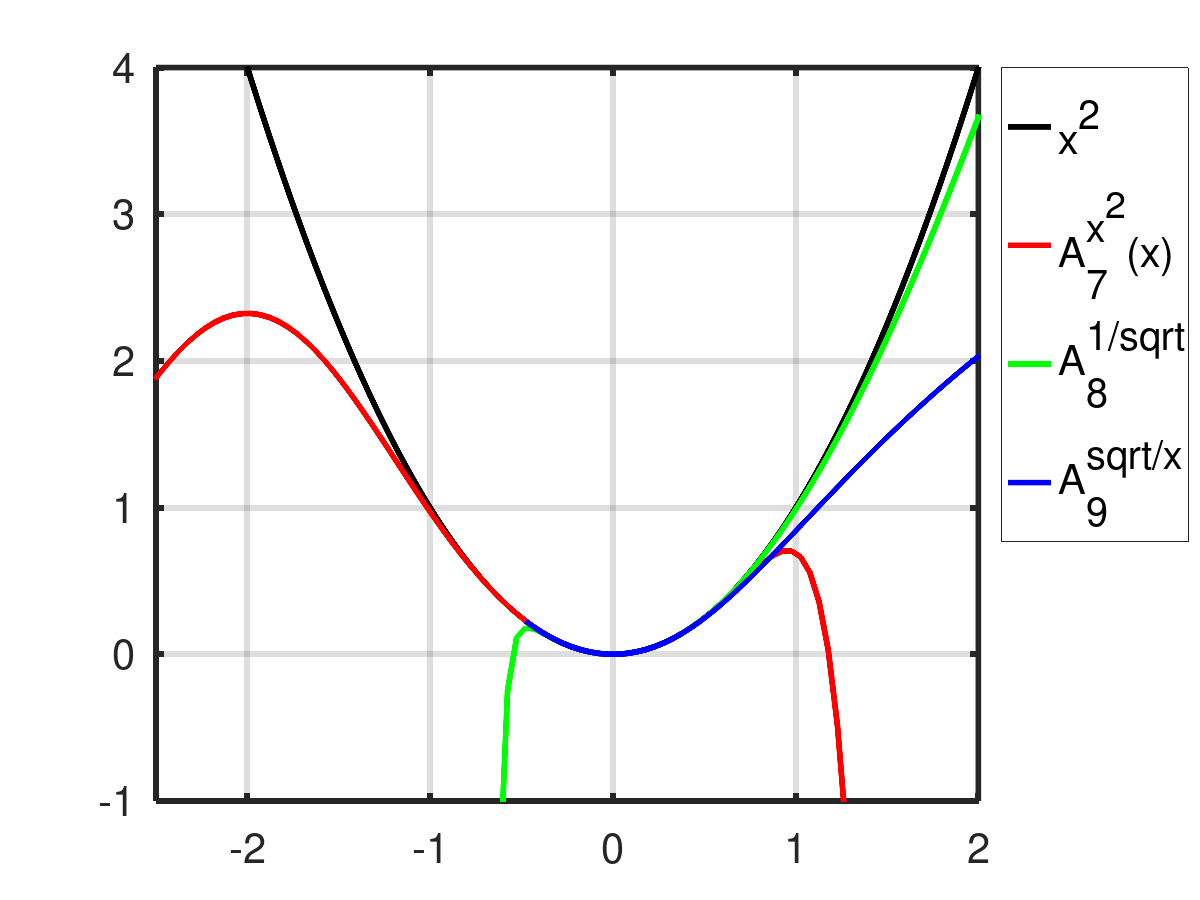}\includegraphics[width=0.5\textwidth]{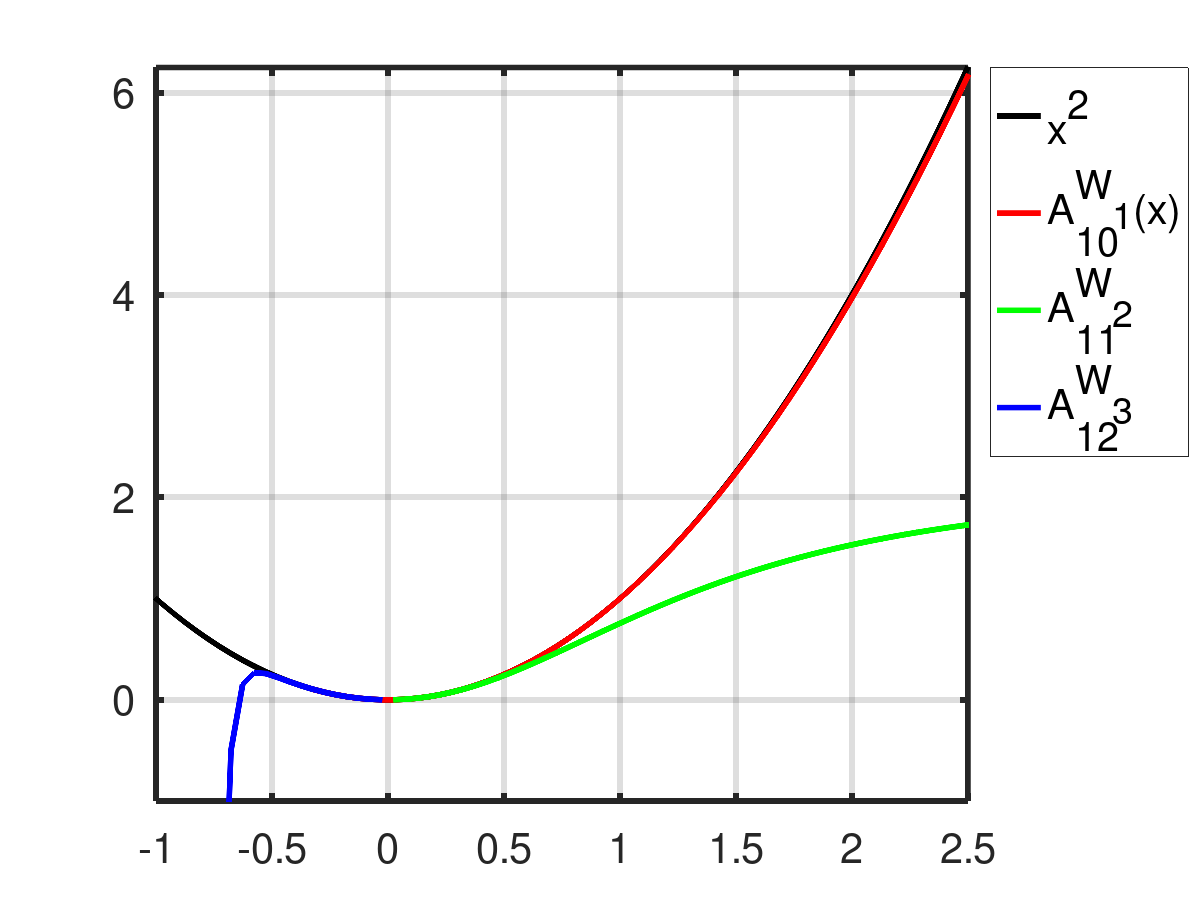}
\par\end{centering}
\begin{centering}
\includegraphics[width=0.5\textwidth]{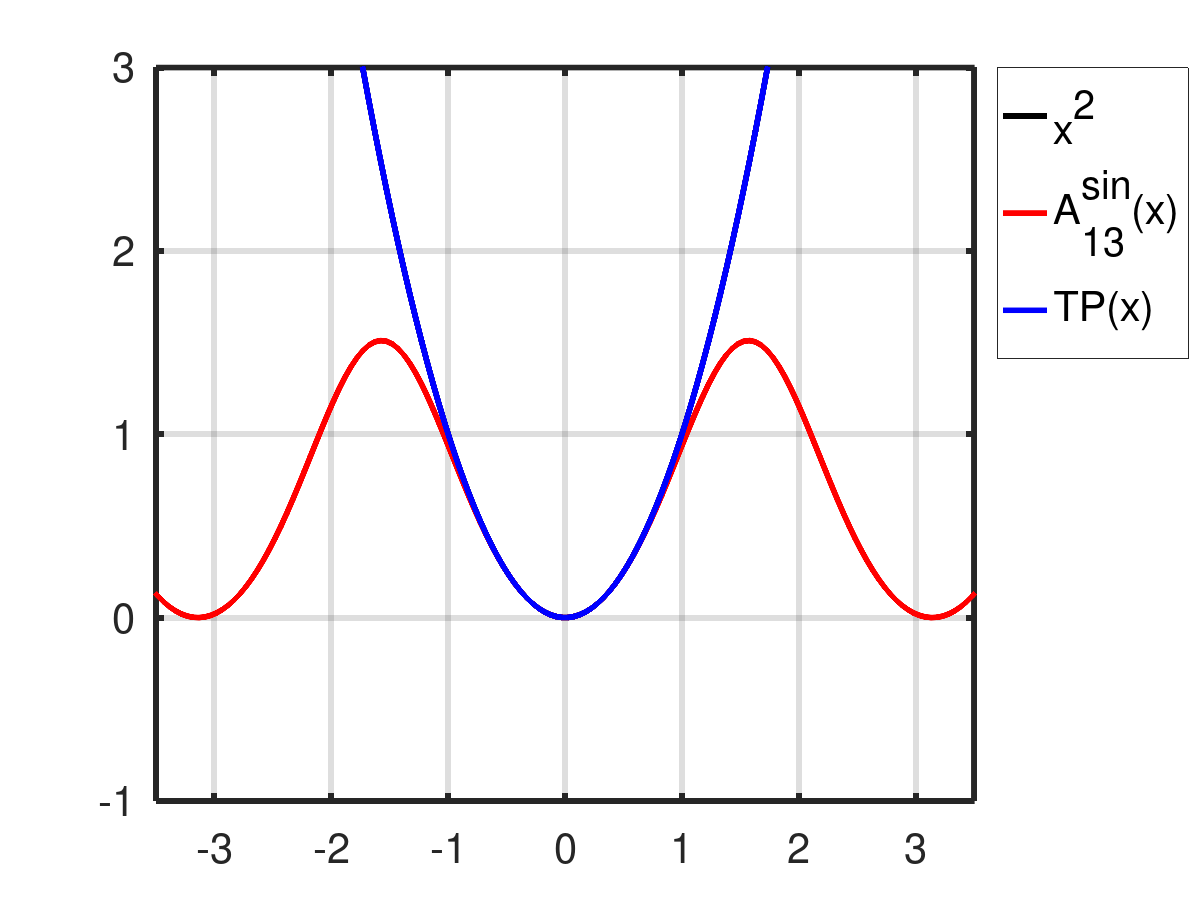}
\par\end{centering}
\caption{The $x^{2}$ function approximated by series (\ref{eq:mainSeries})
with 8 terms and with various $g$ from Sec. \ref{subsec:ExplicitG}. }
\label{FigX2}
\end{figure}

\begin{figure}
\begin{centering}
\includegraphics[width=0.5\textwidth]{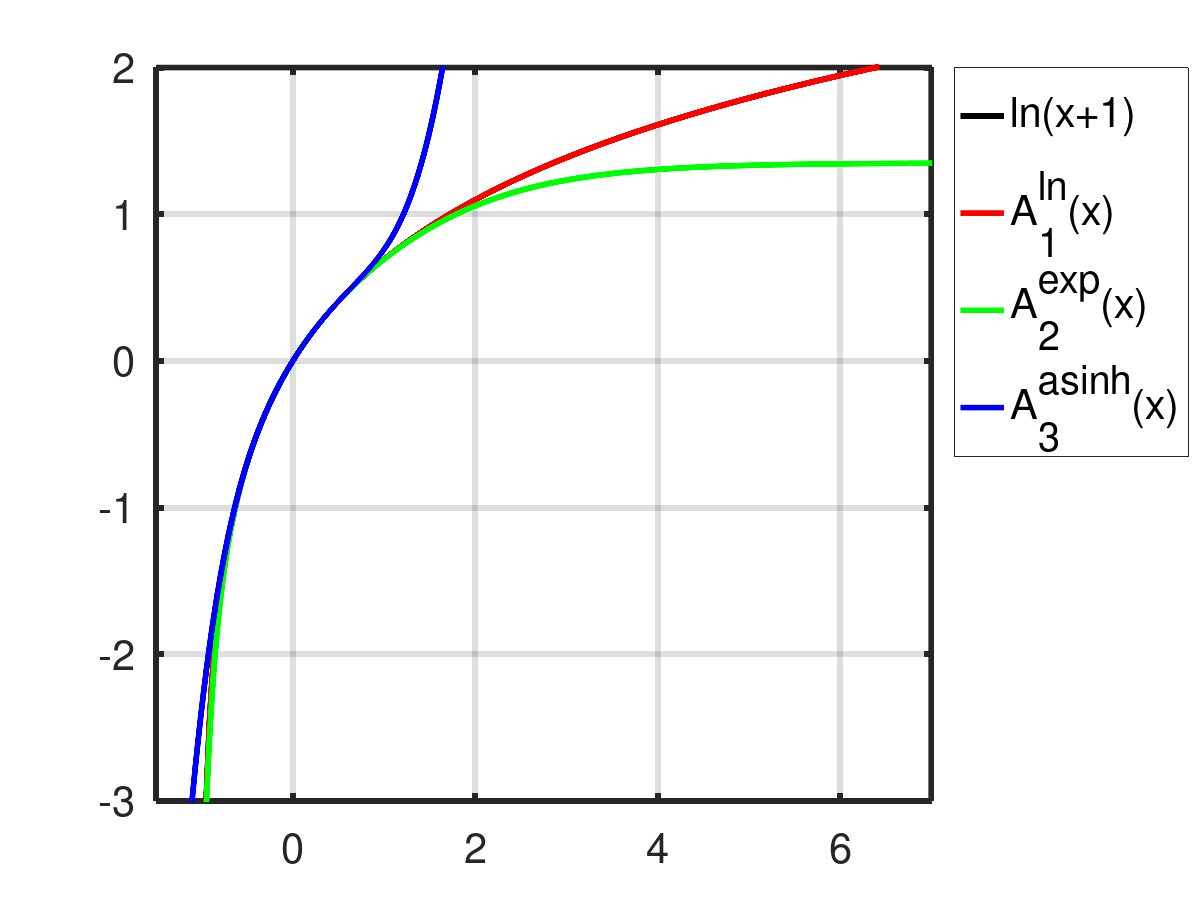}\includegraphics[width=0.5\textwidth]{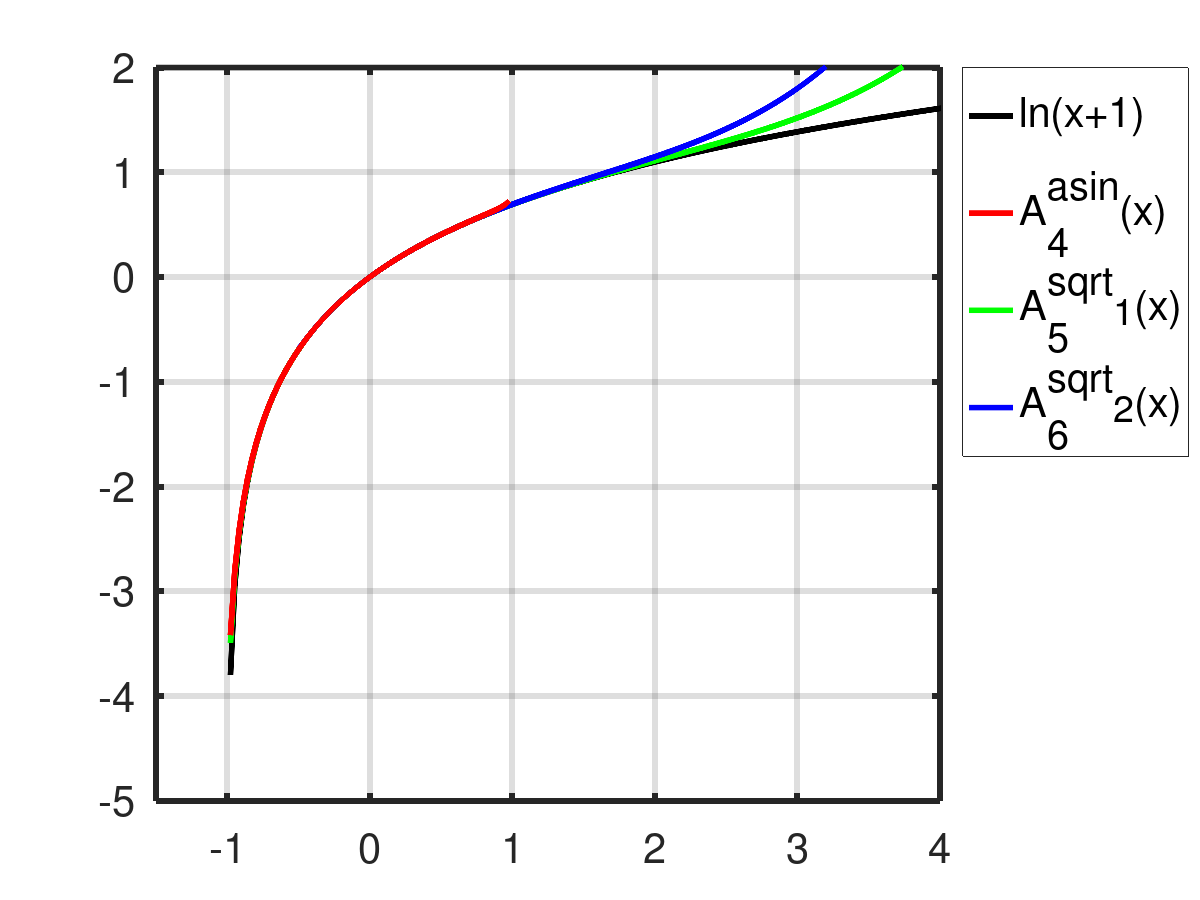}
\par\end{centering}
\begin{centering}
\includegraphics[width=0.5\textwidth]{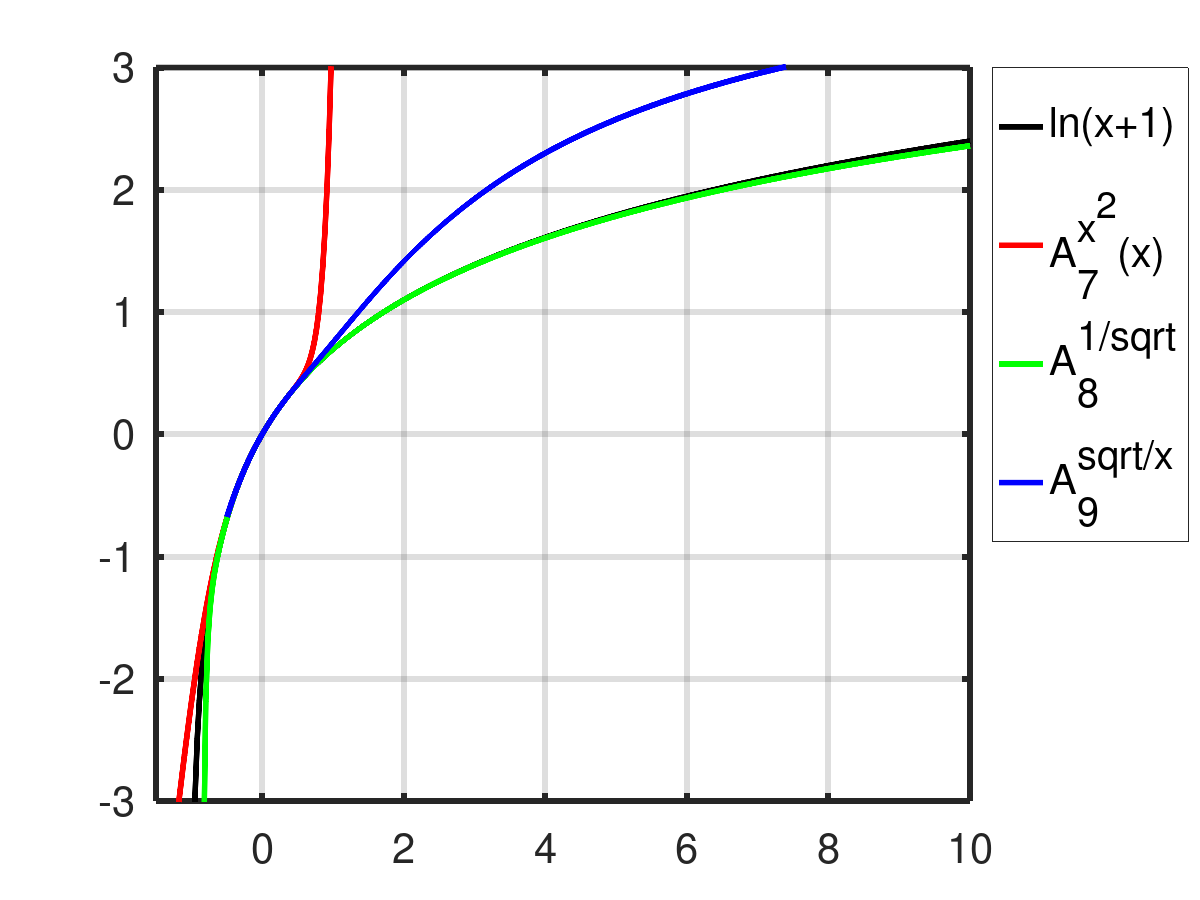}\includegraphics[width=0.5\textwidth]{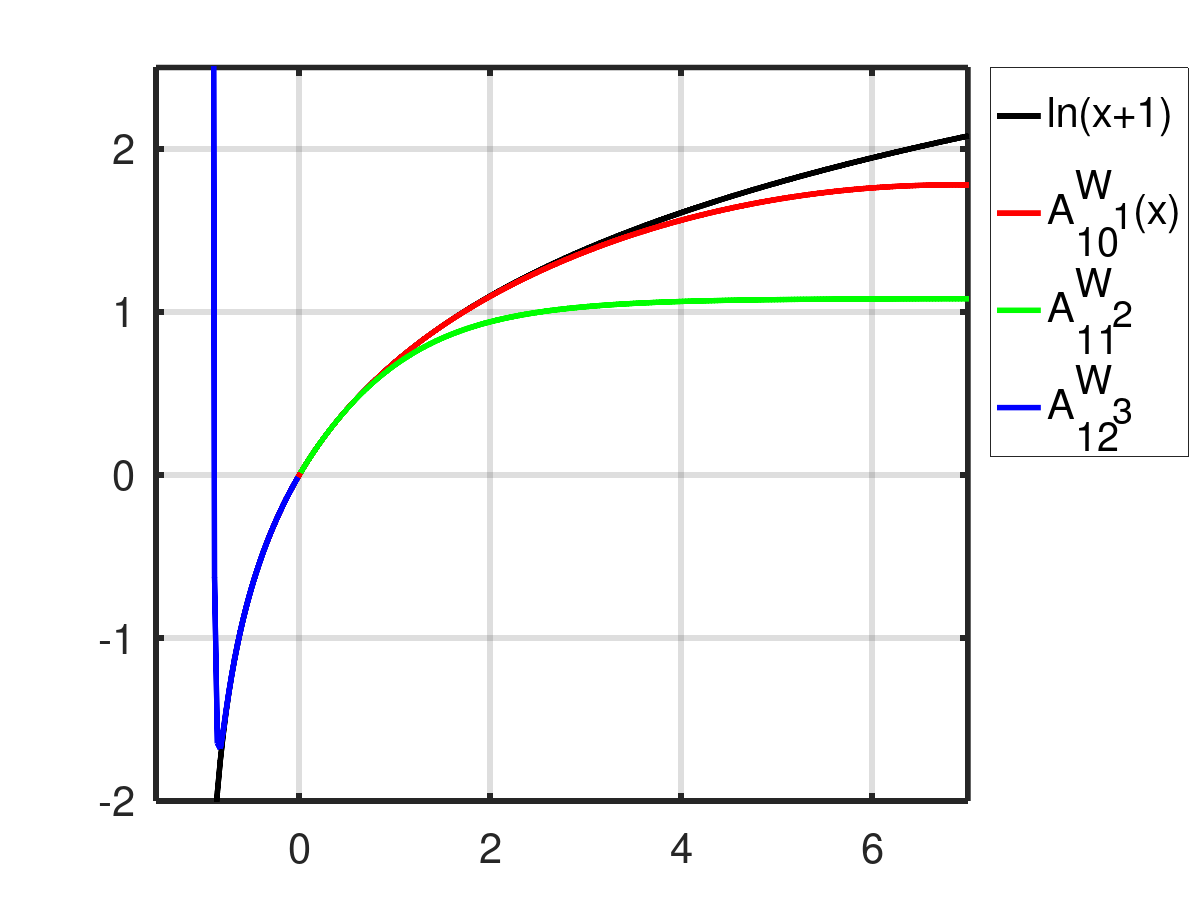}
\par\end{centering}
\begin{centering}
\includegraphics[width=0.5\textwidth]{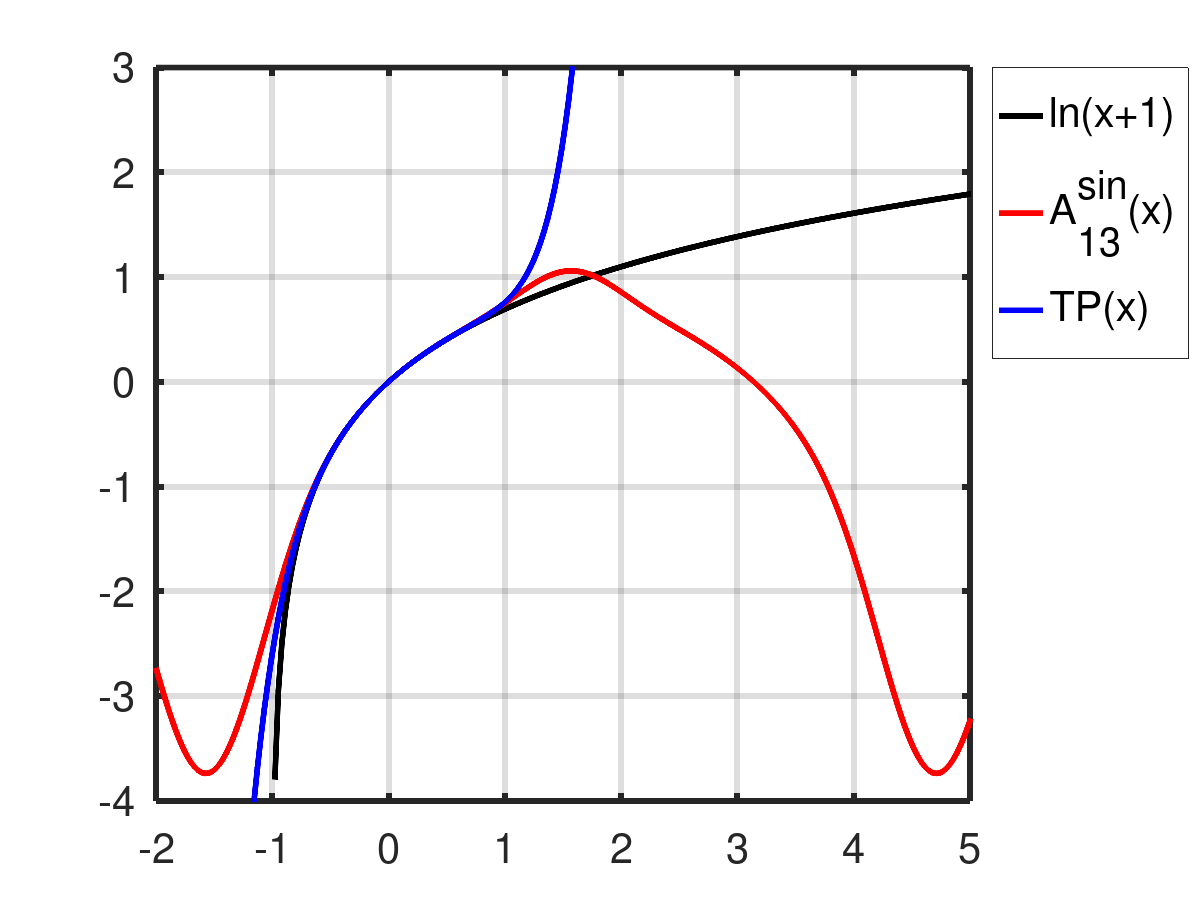}
\par\end{centering}
\caption{The $\ln\left(x+1\right)$ function approximated by series (\ref{eq:mainSeries})
with 8 terms and with various $g$ from Sec. \ref{subsec:ExplicitG}.}
\label{FigLn}
\end{figure}

In Figs. (\ref{FigExp})-(\ref{FigLn}), situated at the end of this
text, we provide plots where four elementary functions $\exp\left(x\right)$,
$\sin\left(x\right)$, $x^{2}$ and $\ln\left(x+1\right)$ are approximated
with expansions based on Eqs. (\ref{eq:app1})-(\ref{eq:app13}),
the value and first seven derivatives are matched. For the sake of
comparison we also include the TP. The numbering subscript of approximations
$\text{A}_{\text{no.}}$ in the legend respects the order in which
the $g$ functions are presented in the Sec. \ref{subsec:ExplicitG}
and the superscript attempts to mimic the function form of $g$ so
as to remind the reader about it. The parametric expressions (\ref{eq:app5}),(\ref{eq:app6}),(\ref{eq:app7})
and (\ref{eq:app10}) are show with parameters $\alpha=2$, $w=1$,
$\left(\alpha=4,\beta=3\right)$ and $w=1$, respectively. Some lines
in the graphs are overlaid, the reason is mostly the fact that the
approximation is exact\footnote{Sin(x) is exactly approximated by (\ref{eq:app13}), $x^{2}$ by (\ref{eq:app5}),(\ref{eq:app6})
and the TP and $\ln\left(x+1\right)$ by (\ref{eq:app1}).}.

\subsection*{Convergence}

Convergence properties can be easily addressed since the substitution
as expressed by the Eq. (\ref{eq:Substitucia}) does not influence
the point-wise behavior. So, considering 
\[
f\left(x\right)=f\left[g^{-1}\left(y\right)\right]\approx\sum_{n=0}^{\infty}a_{n}y^{n},
\]
one applies the standard convergence criteria known from the usual
power series to the coefficient sequence $\{a_{n}\}$ and determines
the radius of convergence $R$ for the variable $y$
\[
\left|y\right|<R\Rightarrow\sum_{n=0}^{\infty}a_{n}y^{n}\text{ converges}.
\]
Then for all $x\in U$, $U=\{x\in\mathbb{R}:\left|g\left(x\right)\right|<R\}$,
the series $\sum_{n=0}^{\infty}a_{n}\left[g\left(x\right)\right]^{n}$converges.

The convergence to the approximated function can also be treated in
this way, for simplicity we assume that we work on an interval $\mathcal{I}$
containing zero where $g$ can be inverted. Writing an equality which
includes the reminder term 
\[
f\left[g^{-1}\left(y\right)\right]=\sum_{n=0}^{M}a_{n}y^{n}+R_{M}\left(y\right),
\]
one can apply the standard criteria known from the Taylor series to
see whether, in a point-wise way, the reminder vanishes with $M\rightarrow\infty$
at some $y_{0}$. If $W$ is the set of all points such that 
\[
y\in W\Rightarrow\sum_{n=0}^{\infty}a_{n}y^{n}=f\left[g^{-1}\left(y\right)\right],
\]
then for all $x\in\mathbb{\mathcal{I}}$ such that $g\left(x\right)\in W$
one has $f\left(x\right)=\sum_{n=0}^{\infty}a_{n}\left[g\left(x\right)\right]^{n}.$

The most difficult part is presumably the application of the standard
criteria to $\{a_{n}\}$, since the expression (\ref{eq:koeficentyVzorec})
is rather complicated (may contain several nested sums).

The convergence criteria can be in a straightforward way extended
to the complex analysis.

\subsection*{Polynomial approximations}

One observes that pure polynomial approximations are in the list:
the parametric expression (\ref{eq:app5}) with $\alpha=1/k,\:k\in\mathbb{N}^{+}$
and the expression (\ref{eq:app7}). It is interesting to realize,
that these expansions in general do not exactly approximate polynomials
with the same number of terms. Since the polynomial coefficients are
in the one-to-one correspondence with the derivatives $d_{k}^{f}$,
the two approximations contain the TP as their lower terms up to $x^{N}$.
In addition, they also contain higher order terms which imply the
deviations from the approximated function if the latter is a polynomial
of the degree $N$.

Further, expansions (\ref{eq:app5}) and (\ref{eq:app6}) with shifted
arguments\footnote{Meaning that the derivatives are evaluated at $x_{0}=1$ and $x_{0}=\omega^{2}/2$,
respectively.} $\mathcal{A}_{5}^{f}\left(x-1\right)$ and $\mathcal{A}_{6}^{f}\left(x-\omega^{2}/2\right)$
allow to build expressions where the integer and/or fractional powers
of $x$ appear. They represent fractional order polynomials which
have already been introduced in the literature and in a special case
are written as
\begin{equation}
F\left(x\right)=\sum_{n=0}^{N}c_{k}\left(x^{\alpha}\right)^{n},\quad\alpha>0,\label{eq:partDegPoly}
\end{equation}
see e.g. Eq. (4) in \cite{fracPoly1} or Eq. (11) in \cite{fracPoly2}.

\subsection*{Applications}

The applications may result from better approximation properties than
what is provided by the TPs. This however depends on the approximated
function, yet some claims are evident, e.g. there are cases where
an approximation proposed here converges beyond the radius of the
convergence of the Taylor series. Indeed, the function $\ln\left(x+1\right)$,
when expanded at zero, can be approximated by the TPs on the interval
$\left(-1,1\right)$ only. By (\ref{eq:app1}) it is approximated
on the whole definition interval exactly and with one term.

To be more fair, we compare the expansion in powers of $g$ from Eq.
(\ref{eq:app8}) with the TP inside its radius of convergence, i.e.
we numerically investigate the approximation of $\ln\left(x+1\right)$
at $x=0.5$. We define $\varDelta_{f}=\left|\ln\left(1.5\right)-f(1.5)\right|$
and we get ($N$ is the number of terms in the series, see (\ref{eq:mainSeries}))
\begin{center}
\begin{tabular}{ccccc}
\hline 
$N$ & 3 & 7 & 10 & 20\tabularnewline
$\varDelta_{\mathcal{A}_{8}^{f}}\approx$ & $6.65\times10^{-4}$ & $3.84\times10^{-7}$ & $1.74\times10^{-9}$ & $3.33\times10^{-16}$\tabularnewline
$\varDelta_{TP}\approx$ & $1.12\times10^{-2}$ & $3.38\times10^{-4}$ & $3.05\times10^{-5}$ & $1.53\times10^{-8}$\tabularnewline
\hline 
\end{tabular}
\par\end{center}

The first few cutoff series for both cases indicate a significant
difference in the rate of convergence in favor of the expansion $\mathcal{A}_{8}^{f}$.

An important disadvantage for an eventual implementation of the series
(\ref{eq:app1})-(\ref{eq:app13}) on a computer might be the time
necessary for computing $g\left(x\right)$ from $x$. To speed up
the evaluation of (\ref{eq:xLikeExpansion}) the Horner's method is
to be used. More importantly, a couple of expansions from Sec. \ref{subsec:ExplicitG}
are based on the square root, which is for several common architectures
implemented as a basic arithmetic operation included into the instruction
set of the processor (often labeled \emph{fsqrt}, see \cite{intel2022}
for x86, \cite{arm2020} for ARM ). This means it can be evaluated
very rapidly which, in combination with possible better convergence
properties, can be a reason for implementing new algorithms to compute
values of some functions.

\begin{figure}
\begin{centering}
\includegraphics[width=0.65\textwidth]{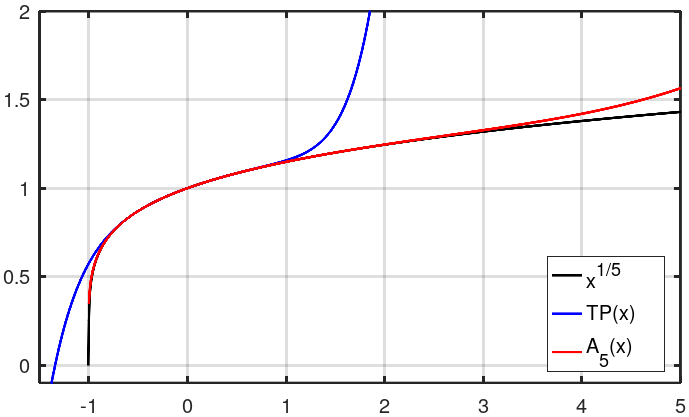}
\par\end{centering}
\caption{The function $\sqrt[5]{x+1}$ approximated by the series $\mathcal{A}_{5}$
built from the powers of $\left[\sqrt{x+1}-1\right]$ and by the TP,
in both cases with 8 terms.}
\label{Fig:5root}

\end{figure}

In this spirit, one potentially interesting application is the computation
of the $m$th root, which is (usually) not a basic instruction of
a processor. Our preliminary tests indicate that $\sqrt[m]{x+1}$
can be for $-1<x$, $1<m\in\mathbb{R}$ computed by the series based
on the expansion in powers of (\ref{eq:app5}), $\alpha=2$,
\begin{itemize}
\item with a significantly higher rate of convergence than have the Taylor
series (within its convergence domain) and
\item on a significantly larger interval (i.e. beyond its convergence domain).
\end{itemize}
Such behavior was observed for all $1<M$ we tested. An example with
the fifth root is shown in Fig. \ref{Fig:5root}. Further investigations
need to be done to confirm our claims on a more rigorous basis.

\section{Summary, conclusion, outlook}

In this text we presented a number of presumably new expansions built
as powers series constructed from functions, we addressed and clarified
the question of their point-wise convergence and mentioned some advantages
they may have in comparison with the Taylor polynomials. These advantages
can represent the reason for their application potential in numerical
evaluation of some functions, the issue however requires more detailed
investigation in the future.

\bibliographystyle{gost2008}
\bibliography{gPowers}
 

\end{document}